\documentclass{amsart}
\usepackage[latin1]{inputenc}
\newtheorem{thm}{Theorem}

\newtheorem{lem}[thm]{Lemma}

\theoremstyle{definition}

 \numberwithin{equation}{section}
\newcommand{\To}{\longrightarrow}
\begin{document}

\thanks{Author supported by FPU grant of MEC of Spain.}
%\date{\today}
\title[]{The unit ball of the Hilbert space in its weak topology}%
\author{Antonio Avilés}

\address{Departamento de Matemáticas\\ Universidad de Murcia\\ 30100 Espinardo (Murcia)\\ Spain}%
\email{avileslo@um.es}

\begin{abstract}
We show that the unit ball of $\ell_p(\Gamma)$ in its weak topology is a
continuous image of $\sigma_1(\Gamma)^\mathbb{N}$ and we deduce some
combinatorial properties of its lattice of open sets which are not shared
by the balls of other equivalent norms when $\Gamma$ is uncountable.\\
\end{abstract}

\maketitle

For a set $\Gamma$ and a real number $1<p<\infty$, the Banach
space $\ell_p(\Gamma)$ is a reflexive space, hence its unit ball
is compact in the weak topology and in fact, it is homeomorphic to
the following closed
subset of the Tychonoff cube $[-1,1]^\Gamma$:\\

$$B(\Gamma) = \left\{ x\in [-1,1]^\Gamma : \sum_{\gamma\in\Gamma}|x_\gamma|\leq
1\right\}.$$

The homeomorphism $h:B_{\ell_p(\Gamma)}\To B(\Gamma)$ is given by
$h(x)_\gamma = sign(x_\gamma)\cdot |x_\gamma|^p$. The spaces
homeomorphic to closed subsets of some $B(\Gamma)$ constitute the
class of uniform Eberlein compacta, introduced by Benyamini and
Starbird~\cite{BenSta}. The space $\sigma_k(\Gamma)$, the compact
subset of $\{0,1\}^\Gamma$ which consists of the functions with at
most $k$ nonzero coordinates ($k$ a positive integer) is an
example of a uniform Eberlein compact. In fact, the following
result was proven in~\cite{BenRudWag}:\\

\begin{thm}[Benyamini, Rudin, Wage]\label{imagsubespuec}
Every uniform Eberlein compact of weight $\kappa$ is a continuous
image of a closed
subset of $\sigma_1(\kappa)^\mathbb{N}$.\\
\end{thm}

In the same paper~\cite{BenRudWag}, it was posed the problem
whether in fact, it was possible to get any uniform Eberlein
compact as a continuous image of the full
$\sigma_1(\Gamma)^\mathbb{N}$. This question was answered in the
negative by Bell \cite{BellRamsey}, by considering the
following property:\\

A compact space $K$ verifies property (Q) if for every uncountable
regular cardinal $\lambda$ and every family
$\{U_\alpha,V_\alpha\}_{\alpha<\lambda}$ of open subsets of $K$
with $\overline{U_\alpha}\subset V_\alpha$ one of the following
two alternatives must hold:
\begin{enumerate}
\item either there exists a set $A\subset\lambda$ with
$|A|=\lambda$ such that $U_\alpha\cap U_\beta = \emptyset$ for
every two different elements $\alpha$ and $\beta$ in $A$, \item or
either there exists a set $A\subset\lambda$ with $|A|=\lambda$
such that $V_\alpha\cap V_\beta \neq \emptyset$ for every two
different
elements $\alpha$ and $\beta$ in $A$.\\
\end{enumerate}

Bell proved in \cite{BellRamsey} that property (Q) is satisfied by
all polyadic spaces, that is, continuous images of
$\sigma_1(\Gamma)^\Lambda$ for any sets $\Gamma$ and $\Lambda$,
(this concept was introduced in \cite{Mrowkapolyadic} and studied
earlier by Gerlits \cite{Gerlitspolyadic}), but he constructed a
uniform Eberlein compact without property (Q). Later, Bell
\cite{BellpropertyB} provided another example of a uniform
Eberlein compact which is not a continuous image of any
$\sigma_1(\Gamma)^\mathbb{N}$ but which is nevertheless polyadic.
Our main result is the following:

\begin{thm}\label{la bola es imagen}
 $B(\Gamma)$ is a continuous image of $\sigma_1(\Gamma)^\mathbb{N}$.\\
\end{thm}

As a consequence, $B(\Gamma)$ satisfies property (Q) as well as
other properties of the same type introduced by Bell in
\cite{BellpropertyB} and \cite{Bellscadic}. However, if $\Gamma$
is uncountable, we show in Theorem \ref{normaequiv} that a
modification of one of the examples of Bell provides an equivalent
norm on $\ell_p(\Gamma)$ whose unit ball is not a continuous image
of $\sigma_1(\Gamma)^\mathbb{N}$, indeed not satisfying property
(Q). In particular, we are showing the existence of equivalent
norms in the nonseparable $\ell_p(\Gamma)$ whose closed unit balls
are not homeomorphic in the weak topology. This contrasts with the
separable case, since the balls of all separable reflexive Banach
spaces are weakly homeomorphic \cite[Theorem 1.1]{Banakhweakball}.
We refer to \cite{Banakhweakball} for information about the
problem whether the balls of
equivalent norms in a Banach space are weakly homeomorphic in the separable case.\\

Proof of Theorem~\ref{la bola es imagen}: For a set $\Delta$ we
will use the notation $B^+(\Delta)=B(\Delta)\cap [0,1]^\Delta$.
 First, we point out that
$B(\Gamma)$ is a continuous image of $B^+(\Gamma)$. Indeed, if we
consider $\Gamma^o = \Gamma\times\{a,b\}$, we have a continuous
surjection $\psi: B^+(\Gamma^o)\To B(\Gamma)$ given by
$\psi(x)_\gamma =
x_{(\gamma,a)}-x_{(\gamma,b)}$.\\

 In a second step, we
apply the standard procedure to express the space $B^+(\Gamma)$ as
a continuous image of a totally disconnected compact $L_0$. We fix
a sequence $(r_n)_{n=0}^\infty$ of positive real numbers such that
$\sum_{n=0}^\infty r_n = 1$ and such that the continuous map
$\phi:\{0,1\}^\mathbb{N}\To [0,1]$ given by $\phi(x) =
\sum_{n=0}^\infty r_n x_n$ is surjective, for example $r_n =
\frac{1}{2^{n+1}}$. We consider the power
$\phi^\Gamma:\{0,1\}^{\Gamma\times\mathbb{N}}\To [0,1]^\Gamma$ and
then we set:
 \begin{eqnarray*} L_0 &=&
(\phi^{\Gamma})^{-1}(B^+(\Gamma)),\\
f &=& \phi^\Gamma|_{L_0},
\end{eqnarray*}
so that $f:L_0\To B^+(\Gamma)$ is a continuous surjection. It will
be convenient to have an explicit description of $L_0$. For
$x\in\{0,1\}^{\Gamma\times\mathbb{N}}$ and $n\in\mathbb{N}$, we
define $N_n(x)=|\{\gamma\in\Gamma : x_{(\gamma,n)}=1\}|$.

\begin{eqnarray*} x\in L_0 &\iff& \phi^\Gamma(x)\in B^+(\Gamma)\\
&\iff& \sum_{\gamma\in\Gamma}\phi^\Gamma(x)_\gamma \leq 1\\
&\iff& \sum_{\gamma\in\Gamma}\sum_{n=0}^\infty r_n x_{(\gamma,n)}\leq 1\\
&\iff& \sum_{n=0}^\infty r_n N_n(x)\leq 1.\\
\end{eqnarray*}

The compact space $L_0$ can be alternatively described as follows. Let $Z$
be a compact subset of $\mathbb{N}^\mathbb{N}$ such that if $\sigma\in Z$
and $\tau_n\leq\sigma_n$ for all $n\in\mathbb{N}$, then $\tau\in Z$.
Associated to such a set $Z$ we construct the following space:

$$\mathcal{K}(Z,\Gamma) = \{x\in\{0,1\}^{\Gamma\times\mathbb{N}} :
(N_n(x))_{n\in\mathbb{N}}\in Z\}.$$

We have that $L_0 =\mathcal{K}(Z_0,\Gamma)$ where $Z_0 =
\left\{s\in\mathbb{N}^\mathbb{N} :\sum_{i\in\mathbb{N}} r_i
s_i\leq 1\right\}$. Note that $Z_0$ is indeed compact since it is
a closed subset of $\prod_{n\in\mathbb{N}}\{0,\ldots,M_n\}$ where
$M_n$ is the integer part of $\frac{1}{r_n}$. The proof will be
complete after the following
lemma:\\

\begin{lem}
Let $Z$ be a compact subset of $\mathbb{N}^\mathbb{N}$ such that
if $\sigma\in Z$ and $\tau_n\leq\sigma_n$ for all
$n\in\mathbb{N}$, then $\tau\in Z$. Then $\mathcal{K}(Z,\Gamma)$
is a continuous image of $\sigma_1(\Gamma)^\mathbb{N}$.\\
\end{lem}

PROOF: First we check that $\mathcal{K}(Z,\Gamma)$ is a closed
subset of $\{0,1\}^{\Gamma\times\mathbb{N}}$ and hence compact.
Namely, if $x\in
\{0,1\}^{\Gamma\times\mathbb{N}}\setminus\mathcal{K}(Z,\Gamma)$,
then $(N_n(x))_{n\in\mathbb{N}}\not\in Z$ and since $Z$ is closed
in $\mathbb{N}^\mathbb{N}$, there is a finite set
$F\subset\mathbb{N}$ such that $\sigma\not\in Z$ whenever
$\sigma_n = N_n(x)$ for all $n\in F$. Indeed, by the definition of
$Z$, if $\tau\in\mathbb{N}^\mathbb{N}$ and $\tau_n\geq\sigma_n$ of
all $n\in F$, also $\tau\not\in Z$. In this case,
$$W =\{y\in\{0,1\}^{\Gamma\times\mathbb{N}} : y_{\gamma,n}=1 \text{ whenever } n\in F\text{ and }x_{\gamma,n}=1\}$$
is a neighborhood which separates $x$ from $\mathcal{K}(Z,\Gamma)$
and this finishes the proof that $\mathcal{K}(Z,\Gamma)$ is
closed. Since $Z$ is compact, for every $n\in\mathbb{N}$ there
exists $M_n\in\mathbb{N}$ such that $\sigma_n\leq M_n$ for all
$\sigma\in Z$. We define the following compact space:

$$L_1 = Z \times \prod_{m\in\mathbb{N}}\prod_{i=0}^{M_m}\sigma_i(\Gamma)$$

Note that $L_1$ is a continuous image of
$\sigma_1(\Gamma)^\mathbb{N}$. On the one hand, since $Z$ is a
metrizable compact, it is a continuous image of
$\{0,1\}^\mathbb{N}$ and in particular of
$\sigma_1(\Gamma)^\mathbb{N}$. On the other hand, for any
$i\in\mathbb{N}$, the space $\sigma_i(\Gamma)$ can be viewed as
the family of all subsets of $\Gamma$ of cardinality at most $i$.
In this way, we consider the continuous surjection
$p:\sigma_1(\Gamma)^i\To\sigma_i(\Gamma)$ given by
$p(x_1,\ldots,x_i)=x_1\cup\cdots\cup x_i$. From the existence of
such a function follows the fact that any countable product of
spaces $\sigma_i(\Gamma)$ is a continuous image of
$\sigma_1(\Gamma)^\mathbb{N}$, and in particular, the second
factor in the expression of $L_1$ is such an image.\\

It remains to define a continuous surjection
$g:L_1\To\mathcal{K}(Z,\Gamma)$. We first fix some notation. An
element of $L_1$ will be written as $(z,x)$ where $z\in Z$ and
$x\in\prod_{m\in\mathbb{N}}\prod_{i=0}^{M_m}\sigma_i(\Gamma)$. At
the same time, such an $x$ is of the form $(x^m)_{m\in\mathbb{N}}$
with $x^m\in \prod_{i=0}^{M_m}\sigma_i(\Gamma)$ and again each
$x^m$ is $(x^{m,i})_{i=1}^{M_m}$ where
$x^{m,i}\in\sigma_i(\Gamma)$. Finally $x^{m,i} =
(x^{m,i}_\gamma)_{\gamma\in\Gamma}\in\sigma_i(\Gamma)\subset\{0,1\}^\Gamma$.
The function
$g:L_1\To\mathcal{K}(Z,\Gamma)\subset\{0,1\}^{\Gamma\times\mathbb{N}}$
is defined as follows:

$$g(z,x)_{\gamma,m} = x^{m,z(m)}_\gamma$$

Observe that $g(x,z)$ maps indeed $L_1$ onto
$\mathcal{K}(Z,\Gamma)$ because for every $m$,
$(x_\gamma^{m,z(m)})_{\gamma\in\Gamma}$ is an arbitrary element of
$\sigma_{z(m)}(\Gamma)$. $\qed$\\

\begin{thm}\label{normaequiv}
Let $\Gamma$ be an uncountable set and $1<p<\infty$. There exists
an equivalent norm on $\ell_p(\Gamma)$ whose unit ball does not
satisfy property (Q) and hence it is not polyadic.\\
\end{thm}

PROOF: This is a variation of an example of
Bell~\cite{BellRamsey}, originally a scattered compact, so that to
make it absolutely convex. We will consider $\omega_1$ as a subset
of $\Gamma$. Let $\phi:\omega_1\To\mathbb{R}$ be a one-to-one map
and
$$G=\{(\alpha,\beta)\in \omega_1\times\omega_1 : \phi(\alpha)<\phi(\beta) \iff \alpha\preceq
\beta\}.$$ We define an equivalent norm on
$\ell_p(\Gamma)\times\ell_p(\Gamma)\sim\ell_p(\Gamma)$ by

$$ \|(x,y)\|' = \sup \{\|x\|_p,\|y\|_p, |x_\alpha|+|y_\beta| :
(\alpha,\beta)\in G\}.$$ and let $K$ be its unit ball considered in its
weak topology. Fix numbers $1<\xi_1<\xi_2<2^{1-\frac{1}{p}}$. The families
of open sets
$$U_\alpha = \{ (x,y)\in K : |x_\alpha|+|y_\alpha|>\xi_2\},\ \ \alpha<\omega_1$$

$$V_\alpha = \{ (x,y)\in
K : |x_\alpha|+|y_\alpha|>\xi_1\},\ \ \alpha<\omega_1$$

 verify that $\overline{U_\alpha}\subset V_\alpha$ and that for any $\alpha,\beta<\omega_1$,
 $U_\alpha\cap U_\beta = \emptyset$ if and only if $(\alpha,\beta)\in G$ if and only if
 $V_\alpha\cap V_\beta = \emptyset$. Namely, if there is some $(x,y)\in V_\alpha\cap
 V_\beta$, then
 $$|x_\alpha|+|y_\alpha|+|x_\beta|+|y_\beta|>\xi_1+\xi_1>2$$
and therefore either $|x_\alpha|+|x_\beta|>1$ or $|y_\alpha|+|y_\beta|>1$
and this implies that $(\alpha,\beta)\not\in G$ since $(x,y)\in K$. On the
other hand, if $(\alpha,\beta)\not\in G$ then the element $(x,y)\in
\ell_p(\Gamma)\times\ell_p(\Gamma)$ which has all coordinates zero except
$x_\alpha = x_\beta=y_\alpha=y_\beta=2^{-\frac{1}{p}}$ lies in
$U_\alpha\cap U_\beta$. Since there is no uncountable well ordered (or
inversely well ordered) subset of $\mathbb{R}$ there is no uncountable
subset $A$ of $\omega_1$ such that $A\times A\subset G$ or $(A\times
A)\cap G = \emptyset$. Therefore, the families $\{U_\alpha\}$ and
$\{V_\alpha\}$ witness the fact that $K$
does not have property (Q) and hence, it is not polyadic.$\qed$\\

The present work was written during a visit to the University of
Warsaw. The author wishes to thank their hospitality, specially to
Witold Marciszewski and Roman Pol, and to Rafa\l\ G\'orak, from
the Polish Academy of Sciences. This work owes very much to the
discussion with them and their suggestions.

\bibliographystyle{amsplain}

\begin{thebibliography}{1}

\bibitem{Banakhweakball}
T.~Banakh, \emph{The topological classification of weak unit balls
of {B}anach
  spaces}, Dissertationes Math. (Rozprawy Mat.) \textbf{387} (2000), 7--35.

\bibitem{BellRamsey}
M.~Bell, \emph{A {R}amsey theorem for polyadic spaces}, Fund.
Math.
  \textbf{150} (1996), no.~2, 189--195.

\bibitem{Bellscadic}
\bysame, \emph{On character and chain conditions in images of
products}, Fund.
  Math. \textbf{158} (1998), no.~1, 41--49.

\bibitem{BellpropertyB}
\bysame, \emph{Polyadic spaces of countable tightness}, Topology
Appl.
  \textbf{123} (2002), no.~3, 401--407.

\bibitem{BenRudWag}
Y.~Benyamini, M.~E. Rudin, and M.~Wage, \emph{Continuous images of
weakly
  compact subsets of {B}anach spaces}, Pacific J. Math. \textbf{70} (1977),
  no.~2, 309--324.

\bibitem{BenSta}
Y.~Benyamini and T.~Starbird, \emph{Embedding weakly compact sets
into
  {H}ilbert space}, Israel J. Math. \textbf{23} (1976), no.~2, 137--141.

\bibitem{Gerlitspolyadic}
J.~Gerlits, \emph{On a generalization of dyadicity}, Studia Sci.
Math. Hungar.
  \textbf{13} (1978), no.~1-2, 1--17 (1981).

\bibitem{Mrowkapolyadic}
S.~Mr{\'o}wka, \emph{Mazur theorem and {$m$}-adic spaces}, Bull.
Acad. Polon.
  Sci. S\'er. Sci. Math. Astronom. Phys. \textbf{18} (1970), 299--305.

\end{thebibliography}

\end{document}